\documentclass[11pt]{article}
\usepackage{amssymb,latexsym,float,geometry,cite,enumerate,setspace}
\newtheorem{theorem}{Theorem}
\newtheorem{corollary}[theorem]{Corollary}
\newtheorem{observation}[theorem]{Observation}

\usepackage{lineno}

\begin{document}

\onehalfspace

\title{Equality of Distance Packing Numbers}

\author{Felix Joos and Dieter Rautenbach}

\date{}

\maketitle

\vspace{-1cm}

\begin{center}
Institut f\"{u}r Optimierung und Operations Research, 
Universit\"{a}t Ulm, Ulm, Germany\\
\{\texttt{felix.joos, dieter.rautenbach}\}\texttt{@uni-ulm.de}
\end{center}

\begin{abstract}
We characterize the graphs for which the independence number equals the packing number.
As a consequence we obtain simple structural descriptions of the graphs for which 
(i) the distance-$k$-packing number equals the distance-$2k$-packing number,
and
(ii) the distance-$k$-matching number equals the distance-$2k$-matching number.
This last result considerably simplifies and extends previous results of Cameron and Walker 
(The graphs with maximum induced matching and maximum matching the same size, Discrete Math. 299 (2005) 49-55).
For positive integers $k_1$ and $k_2$ with
$k_1<k_2$ and $\lceil(3k_2+1)/2\rceil\leq 2k_1+1$,
we prove that it is NP-hard to determine for a given graph whether its 
distance-$k_1$-packing number equals its distance-$k_2$-packing number.
\end{abstract}


\section{Introduction}

Induced matchings in graphs were introduced by Stockmeyer and Vazirani \cite{stva} as a variant of ordinary matchings.
While the structure and algorithmic properties of ordinary matchings are well understood \cite{lopl},
induced matchings are algorithmically very hard \cite{ca,stva,dadelo}.
Many efficient algorithms for finding maximum induced matchings exploit the fact that induced matchings correspond to independent sets of the square of the line graph \cite{brmo,casrta,ch,ca2,gole}.
In \cite{koro} Kobler and Rotics showed that the graphs where the matching number and the induced matching number coincide can be recognized efficiently. Their result was extended by Cameron and Walker \cite{cawa} who gave a complete structural description of these graph.
In \cite{dujoperaso} we generalized some results from \cite{koro,cawa} to distance-$k$-matchings and simplified the original proofs.
In the present paper we present much more general results systematically exploiting the above-mentioned relation between matchings and independent sets in line graphs. 
Our main result is a very simple characterization of the graphs for which the independence number equals the packing number.
An immediate consequence of this result is a complete structural description of the graphs for which the distance-$k$-matching number equals the distance-$2k$-matching number.
It follows immediately that such graphs can be recognized by a very simple efficient algorithm.
We establish further results relating distance packing numbers and discuss related open problems.

Before we proceed to the results, we recall some terminology.
We consider finite, simple, and undirected graphs.
Let $G$ be a graph.
A set $P$ of vertices of $G$ is a {\it $k$-packing of $G$} for some positive integer $k$
if every two distinct vertices in $P$ have distance more than $k$ in $G$.
The {\it $k$-packing number $\rho_k(G)$ of $G$} is the maximum cardinality of a $k$-packing of $G$,
and a $k$-packing of cardinality $\rho_k(G)$ is {\it maximum}.
Using this terminology, independent sets correspond to $1$-packings and 
the independence number $\alpha(G)$ coincides with $\rho_1(G)$.
We denote the {\it line graph of $G$} by $L(G)$
and the {\it $k$-th power} of $G$ for some positive integer $k$ by $G^k$.
Since matchings of $G$ correspond to independent sets of $L(G)$, 
the matching number $\nu(G)$ equals $\rho_1(L(G))$.
Similarly, 
since induced matchings of $G$ correspond to $2$-packings of $L(G)$, 
the induced matching number $\nu_2(G)$ equals $\rho_2(L(G))$.
More generally, a set $M$ of edges of $G$ is a {\it $k$-matching of $G$}
if it is a $k$-packing of $L(G)$.
The {\it $k$-matching number $\nu_k(G)$} and {\it maximum $k$-matchings} 
are defined in the obvious way.
Clearly, a set $P$ is a $k_1$-packing of $G^{k_2}$
for some positive integers $k_1$ and $k_2$ if and only it is a $k_1k_2$-packing of $G$, that is,
$\rho_{k_1}(G^{k_2})=\rho_1(G^{k_1k_2})$.
A vertex $u$ of $G$ is {\it simplicial} 
if $N_G[u]$ is complete. Two distinct vertices $u$ and $v$ of $G$ are {\it twins} if $N_G[u]=N_G[v]$.
Let $S(G)$ be the set of simplicial vertices of $G$.
Let ${\cal S}(G)$ be the partition of $S(G)$ where two simplicial vertices belong to the same partite set if and only if they are twins.
A {\it transversal of ${\cal S}(G)$} is a set of simplicial vertices that contains exactly one vertex from each partite set of the partition ${\cal S}(G)$.
Note that the subgraph of $G$ induced by $S(G)$ is a union of cliques, and that ${\cal S}(G)$ is the collection of the vertex sets of these cliques. In particular, every transversal of ${\cal S}(G)$ is independent.

\section{Results}

We immediately proceed to the characterization of the graphs for which the independence number equals the packing number.

\begin{theorem}\label{theorem1}
A graph $G$ satisfies $\rho_1(G)=\rho_2(G)$ if and only if 
\begin{enumerate}[(i)]
\item a set of vertices of $G$ is a maximum $2$-packing if and only if it is a transversal of ${\cal S}(G)$, and 
\item for every transversal $P$ of ${\cal S}(G)$, the sets $N_G[u]$ for $u$ in $P$ partition $V(G)$.
\end{enumerate}
\end{theorem}
{\it Proof:} Let $G$ be a graph.

In order to prove the sufficiency, let $G$ satisfy (i) and (ii).
Let $P$ be a transversal of ${\cal S}(G)$. 
By (i), we have $|P|=\rho_2(G)$. 
By (ii) and since $P\subseteq S(G)$, we obtain that $\{ N_G[u]: u\in P\}$ is a partition of $V(G)$ into complete sets. Since every $1$-packing contains at most one vertex from each complete set, this implies $\rho_1(G)\leq |P|$.
Since $\rho_1(G)\geq \rho_2(G)$, it follows $\rho_1(G)=\rho_2(G)$.

In order to prove the necessity, let $G$ satisfy $\rho_1(G)=\rho_2(G)$.
Let $P$ be a maximum $2$-packing.
If some vertex $u$ in $P$ has two non-adjacent neighbors $v$ and $w$, 
then $(P\setminus \{ u\})\cup \{ v,w\}$ is a $1$-packing with more vertices than $P$, 
which is a contradiction.
Hence all vertices in $P$ are simplicial.
Since no two vertices in $P$ are adjacent, 
the set $P$ is contained in some transversal $Q$ of ${\cal S}(G)$. 
Since $Q$ is a $1$-packing,
we obtain $\rho_2(G)=|P|\leq |Q|\leq \rho_1(G)=\rho_2(G)$,
that is, $P=Q$, which implies in particular that $P$ is a transversal of ${\cal S}(G)$.
If $V(G)\setminus \bigcup_{u\in P}N_G[u]$ contains a vertex $v$, then $P\cup \{ v\}$ is $1$-packing with more vertices than $P$, which is a contradiction.
Hence $\{ N_G[u]: u\in P\}$ is a partition of $V(G)$ into complete sets.
Since for every transversal $P'$ of ${\cal S}(G)$, the partition $\{ N_G[u']: u'\in P'\}$ equals the partition $\{ N_G[u]: u\in P\}$,
it follows that every transversal of ${\cal S}(G)$ is a maximum $2$-packing.
Altogether, (i) and (ii) follow. $\Box$

\medskip

\noindent By considering suitable powers of the underlying graph,
we obtain the following.

\begin{corollary}\label{corollary1}
A graph $G$ satisfies $\rho_k(G)=\rho_{2k}(G)$ for some positive integer $k$ if and only if 
\begin{enumerate}[(i)]
\item a set of vertices of $G$ is a maximum $2k$-packing if and only if it is a transversal of ${\cal S}(G^k)$, and 
\item for every transversal $P$ of ${\cal S}(G^k)$, the sets $N_{G^k}[u]$ for $u$ in $P$ partition $V(G)$.
\end{enumerate}
\end{corollary}
By Corollary \ref{corollary1}, it is algorithmically very easy to recognize the graphs $G$ with $\rho_k(G)=\rho_{2k}(G)$.

In view of Theorem \ref{theorem1} and Corollary \ref{corollary1},
it makes sense to consider the equality of distance packing numbers
$\rho_{k_1}(G)$ and $\rho_{k_2}(G)$ where $k_1<k_2$ are positive integers that do not satisfy $k_2=2k_1$.
Our next observation shows that for $k_2>2k_1$ such graphs are not very interesting.

\begin{observation}\label{observation1}
If $k_1$ and $k_2$ are positive integers with $k_2>2k_1$
and $G$ is a connected graph with $\rho_{k_1}(G)=\rho_{k_2}(G)$, 
then $\rho_{k_1}(G)=\rho_{k_2}(G)=1$.
\end{observation}
{\it Proof:} Let $G$ be a graph that satisfies $\rho_{k_1}(G)=\rho_{k_2}(G)$.
Let $P$ be a maximum $k_2$-packing. 
For a contradiction, we assume that $P$ has more than one element.
Let $u$ be a vertex in $P$.
Since $P$ has more than one element, 
there is a vertex $v$ at distance $k_1+1$ from $u$.
Since $k_2+1\geq 2(k_1+1)$, every vertex in $P$ has distance more than $k_1$ from $v$.
Now $P\cup \{ v\}$ is a $k_1$-packing, 
which is a contradiction.
This completes the proof. $\Box$

\medskip

\noindent Now we consider the case $k_1<k_2<2k_1$
and show that already the smallest possible choice, $k_1=2$ and $k_2=3$, leads to graphs that will most likely not have a nice structural description.

\begin{figure}[t]
\begin{center}
\unitlength 1mm 
\linethickness{0.4pt}
\ifx\plotpoint\undefined\newsavebox{\plotpoint}\fi 
\begin{picture}(11,19)(0,0)
\put(0,5){\circle*{2}}
\put(10,5){\circle*{2}}
\put(10,15){\circle*{2}}
\put(0,15){\circle*{2}}
\put(10,15){\line(-1,0){10}}
\put(0,15){\line(0,-1){10}}
\put(0,5){\line(1,0){10}}
\put(10,5){\line(0,1){10}}
{\footnotesize
\put(0,1){\makebox(0,0)[cc]{$x_i$}}
\put(10,1){\makebox(0,0)[cc]{$\bar{x}_i$}}
\put(0,19){\makebox(0,0)[cc]{$\bar{x}_i'$}}
\put(10,19){\makebox(0,0)[cc]{$x_i'$}}
}\end{picture}
\hspace{2cm}
\linethickness{0.4pt}
\ifx\plotpoint\undefined\newsavebox{\plotpoint}\fi 
\begin{picture}(65,30)(0,0)
\put(60,5){\circle*{2}}
\put(40,5){\circle*{2}}
\put(60,15){\circle*{2}}
\put(30,15){\circle*{2}}
\put(50,15){\circle*{2}}
\put(20,15){\circle*{2}}
\put(40,15){\circle*{2}}
\put(10,15){\circle*{2}}
\put(60,25){\circle*{2}}
\put(40,25){\circle*{2}}
\put(50,10){\circle*{2}}
\put(50,20){\circle*{2}}
\put(60,5){\line(-2,1){10}}
\put(50,10){\line(2,1){10}}
\put(60,15){\line(-1,0){10}}
\put(50,15){\line(1,-1){10}}
\put(60,5){\line(-2,3){10}}
\put(50,20){\line(2,-1){10}}
\put(50,20){\line(2,1){10}}
\put(60,25){\line(-1,-1){10}}
\put(50,10){\line(2,3){10}}
\put(40,5){\line(2,1){10}}
\put(50,10){\line(-2,1){10}}
\put(40,15){\line(2,1){10}}
\put(50,20){\line(-2,1){10}}
\put(40,25){\line(1,-1){10}}
\put(50,15){\line(-1,-1){10}}
\put(40,5){\line(1,0){20}}
\put(40,25){\line(1,0){20}}
\qbezier(40,15)(50,11)(60,15)
\put(0,15){\circle*{2}}
\put(0,15){\line(1,0){40}}
\put(40,25){\line(-1,-1){10}}
\put(30,15){\line(1,-1){10}}
{\footnotesize
\put(66,5){\makebox(0,0)[]{$x'(C)$}}
\put(66,25){\makebox(0,0)[]{$z'(C)$}}
\put(66,15){\makebox(0,0)[]{$y'(C)$}}
\put(38,2){\makebox(0,0)[]{$x(C)$}}
\put(38,18){\makebox(0,0)[]{$y(C)$}}
\put(38,28){\makebox(0,0)[]{$z(C)$}}
\put(0,18){\makebox(0,0)[]{$a(C)$}}
\put(28,18){\makebox(0,0)[]{$b(C)$}}
}
\end{picture}
\caption{On the left, the cycle $G(x_i):x_i\bar{x}_ix_i'\bar{x}'_ix_i$ created for the variable $x_i$.
On the right, the graph $G(C)$ created for a clause $C$ with literals $x$, $y$, and $z$, that is, $C=x\vee y\vee z$ and $x,y,z\in \{ x_i:i\in [n]\}\cup \{ \bar{x}_i:i\in [n]\}$.}\label{fig1}
\end{center}
\end{figure}
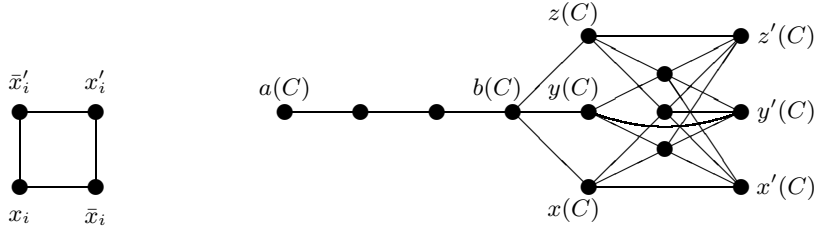

\begin{figure}[t]
\begin{center}
\unitlength 1mm 
\linethickness{0.4pt}
\ifx\plotpoint\undefined\newsavebox{\plotpoint}\fi 
\begin{picture}(91,79)(0,0)
\put(60,30){\circle*{2}}
\put(40,30){\circle*{2}}
\put(60,40){\circle*{2}}
\put(30,40){\circle*{2}}
\put(50,40){\circle*{2}}
\put(20,40){\circle*{2}}
\put(40,40){\circle*{2}}
\put(10,40){\circle*{2}}
\put(0,40){\circle{3}}
\put(40,30){\circle{3}}
\put(60,50){\circle*{2}}
\put(40,50){\circle*{2}}
\put(50,35){\circle*{2}}
\put(50,45){\circle*{2}}
\put(60,30){\line(-2,1){10}}
\put(50,35){\line(2,1){10}}
\put(60,40){\line(-1,0){10}}
\put(50,40){\line(1,-1){10}}
\put(60,30){\line(-2,3){10}}
\put(50,45){\line(2,-1){10}}
\put(50,45){\line(2,1){10}}
\put(60,50){\line(-1,-1){10}}
\put(50,35){\line(2,3){10}}
\put(40,30){\line(2,1){10}}
\put(50,35){\line(-2,1){10}}
\put(40,40){\line(2,1){10}}
\put(50,45){\line(-2,1){10}}
\put(40,50){\line(1,-1){10}}
\put(50,40){\line(-1,-1){10}}
\put(40,30){\line(1,0){20}}
\put(40,50){\line(1,0){20}}
\qbezier(40,40)(50,36)(60,40)
\put(0,40){\line(1,0){40}}
\put(40,50){\line(-1,-1){10}}
\put(30,40){\line(1,-1){10}}
\put(80,5){\circle*{2}}
\put(80,35){\circle*{2}}
\put(80,65){\circle*{2}}
\put(90,5){\circle*{2}}
\put(90,35){\circle*{2}}
\put(90,65){\circle*{2}}
\put(90,15){\circle*{2}}
\put(90,45){\circle*{2}}
\put(90,75){\circle*{2}}
\put(80,15){\circle*{2}}
\put(80,45){\circle*{2}}
\put(80,75){\circle*{2}}
\put(90,15){\line(-1,0){10}}
\put(90,45){\line(-1,0){10}}
\put(90,75){\line(-1,0){10}}
\put(80,15){\line(0,-1){10}}
\put(80,45){\line(0,-1){10}}
\put(80,75){\line(0,-1){10}}
\put(80,5){\line(1,0){10}}
\put(80,35){\line(1,0){10}}
\put(80,65){\line(1,0){10}}
\put(90,5){\line(0,1){10}}
\put(90,35){\line(0,1){10}}
\put(90,65){\line(0,1){10}}
{\footnotesize
\put(80,1){\makebox(0,0)[cc]{$x_1$}}
\put(80,5){\circle{3}}
\put(80,35){\circle{3}}
\put(90,65){\circle{3}}
\put(80,31){\makebox(0,0)[cc]{$x_2$}}
\put(80,61){\makebox(0,0)[cc]{$x_4$}}
\put(90,1){\makebox(0,0)[cc]{$\bar{x}_1$}}
\put(90,31){\makebox(0,0)[cc]{$\bar{x}_2$}}
\put(90,61){\makebox(0,0)[cc]{$\bar{x}_4$}}
\put(77,18){\makebox(0,0)[cc]{$\bar{x}_1'$}}
\put(80,49){\makebox(0,0)[cc]{$\bar{x}_2'$}}
\put(80,79){\makebox(0,0)[cc]{$\bar{x}_4'$}}
\put(90,19){\makebox(0,0)[cc]{$x_1'$}}
\put(90,49){\makebox(0,0)[cc]{$x_2'$}}
\put(90,79){\makebox(0,0)[cc]{$x_4'$}}
\put(0,40){\circle*{2}}
\put(60,25){\makebox(0,0)[cc]{$x_1'(C_1)$}}
\put(57,55){\makebox(0,0)[cc]{$\bar{x}_4'(C_1)$}}
\put(60,30){\line(2,-1){30}}
\put(60,40){\line(4,1){20}}
\put(80,75){\line(-4,-5){20}}
\put(63,36){\makebox(0,0)[cc]{$\bar{x}_2'(C_1)$}}
}
\end{picture}
\caption{The edges between $G(C_1)$ and $G(x_1)\cup G(x_2)\cup G(x_4)$ 
created for the clause $C_1=x_1\vee \bar{x}_2\vee \bar{x}_4$.
If a satisfying truth assignment sets $x_1$ and $x_2$ to true and $x_4$ to false,
the encircled vertices indicate elements of the corresponding $3$-packing.}\label{fig2}
\end{center}
\end{figure}
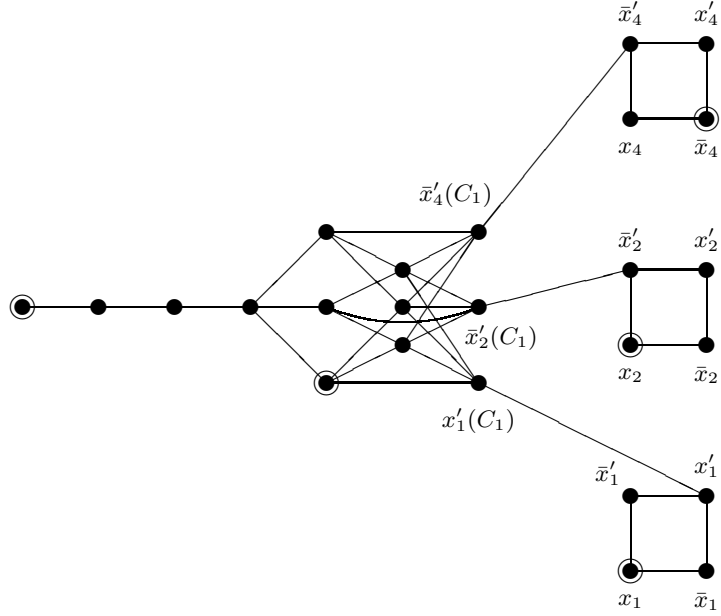

\begin{theorem}\label{theorem2}
It is NP-hard to determine for a given graph $G$ whether $\rho_2(G)=\rho_3(G)$.
\end{theorem}
{\it Proof:} We describe a reduction from 3SAT to the considered problem. Therefore, let $f$ be a 3SAT instance
with $m$ clauses $C_1,\ldots,C_m$ 
over $n$ boolean variables $x_1,\ldots,x_n$.
We construct a graph $G$ whose order is polynomially bounded in terms of $n$ and $m$ 
such that $f$ is satisfiable if and only if $\rho_2(G)=\rho_3(G)$.
For every variable $x_i$, we create a cycle $G(x_i):x_i\bar{x}_ix_i'\bar{x}'_ix_i$ of length $4$
as shown in the left of Figure \ref{fig1}.
For every clause $C_j$, we create a copy $G(C_j)$ of the graph in the right of Figure \ref{fig1} 
and denote its vertices as explained in the caption.
All graphs $G(x_i)$ and $G(C_j)$ created so far are disjoint.
For every clause $C$ with literals $x$, $y$, and $z$, 
we create the three edges
$x'(C)x'$,
$y'(C)y'$, and
$z'(C)z'$.
If, for example, $C_1=x_1\vee \bar{x}_2\vee \bar{x}_4$, 
then these are the edges
$x_1(C)x_1'$,
$\bar{x}_2(C)\bar{x}_2'$, and
$\bar{x}_4(C)\bar{x}_4'$
as shown in Figure \ref{fig2}.
This completes the description of $G$.
It is easy to verify that $\rho_2(G(x_i))=1$ and $\rho_2(G(C_j))=2$,
which implies that $\rho_2(G)\leq n+2m$.
Since 
$\{ a(C_j):j\in [m]\}
\cup \{ b(C_j):j\in [m]\}
\cup \{ x_i:i\in [n]\}$
is a $2$-packing of cardinality $n+2m$,
we obtain $\rho_2(G)=n+2m$.
It remains to prove that $f$ is satisfiable if and only if $\rho_3(G)=n+2m$.

First, we assume that $f$ is satisfiable and consider a satisfying truth assignment.
For every clause $C_j$, we select a true literal $t_j$.
Note that there may be several choices for $t_j$.
Now, by construction, the set
$\{ a(C_j):j\in [m]\}
\cup \{ t_j(C):j\in [m]\}
\cup \{ x_i:i\in [n]\mbox{ and $x_i$ is true}\}
\cup \{ \bar{x}_i:i\in [n]\mbox{ and $x_i$ is false}\}$
is a $3$-packing of cardinality $n+2m$,
which implies $\rho_3(G)=n+2m$.

Next, we assume that $\rho_3(G)=n+2m$.
Let $P$ be a maximum $3$-packing.
Since $\rho_2(G(x_i))=1$ and $\rho_2(G(C_j))=2$,
it follows that 
$P$ contains 
exactly one vertex from each $G(x_i)$
and 
exactly two vertices from each $G(C_j)$.
Clearly, we may assume that for every $i\in [n]$, the set $P$ contains exactly one of the two vertices $x_i$ and $\bar{x}_i$ of the cycle $G(x_i)$.
Similarly, we may assume that for every $j\in [m]$, the set $P$ contains the vertex $a(C_j)$ and exactly one of the three vertices $x(C_j)$, $y(C_j)$, and $z(C_j)$
where $x$, $y$, and $z$ are the three literals in $C_j$.
See Figure \ref{fig2} for an illustration.
We consider the assignment of truth values where the variable $x_i$ is set to true exactly if the vertex $x_i$ belongs to $P$.
If $C$ is a clause and $x\in \{ x_i,\bar{x}_i\}$ is a literal in $C$ such that $P$ contains $x(C)$, then $x'(C)$ is adjacent to the vertex $x'$ of $G(x_i)$, and hence $P$ cannot contain the vertex $x$ of $G(x_i)$.
More specifically, 
if $x=x_i$, then $P$ contains the vertex $x_i$ of $G(x_i)$, 
which means that $x_i$ is set to true,
and 
if $x=\bar{x}_i$, then $P$ contains the vertex $\bar{x}_i$ of $G(x_i)$, 
which means that $x_i$ is set to false.
Altogether, it follows that the truth assignment defined above satisfies $f$.
This completes the proof. $\Box$

\medskip

\noindent A simple modification of the construction in the proof of Theorem \ref{theorem2} allows to establish the following.

\begin{corollary}\label{coro NPhard}
Let $k_1$ and $k_2$ be positive integers with
$k_1<k_2$ and $\lceil(3k_2+1)/2\rceil\leq 2k_1+1$.
It is NP-hard to determine for a given graph $G$
whether $\rho_{k_1}(G)=\rho_{k_2}(G)$.
\end{corollary}
{\it Proof:}
We apply the following modifications
to the graph $G$ 
constructed in the proof of Theorem \ref{theorem2}.
\begin{itemize} 
\item For every variable $x_i$,
subdivide each of the two edges 
$x_i\bar{x}'_i$
and 
$\bar{x}_ix'_i$
exactly $\lceil\frac{k_2}{2}\rceil-2$ times.
\item For each clause $C_j$ with literals $x$, $y$, and $z$,
\begin{itemize} 
\item subdivide the edge incident with $a(C_j)$ 
exactly $k_2-3$ times, and
\item subdivide each of the three edges 
$b(C_j)x(C_j)$,
$b(C_j)y(C_j)$, and
$b(C_j)z(C_j)$
exactly $\lfloor\frac{k_2}{2}\rfloor-1$ times.
\end{itemize} 
\end{itemize} 
Note that after these modifications, 
the distance between $x_i$ and $\bar{x}'_i$
as well as between $\bar{x}_i$ and $x'_i$
is $\lceil\frac{k_2}{2}\rceil-1$,
the distance between $a(C_j)$ and $b(C_j)$ is $k_2$,
and the distance between $b(C_j)$ and each of 
$x(C_j)$, $y(C_j)$, and $z(C_j)$ is 
$\lfloor\frac{k_2}{2}\rfloor$.
Renaming the three neighbors of $b(C_j)$ that do not lie on the path to $a(C_j)$ as $x(C_j)$, $y(C_j)$, and $z(C_j)$,
and repeating the very same argument as in the proof of Theorem \ref{theorem2}, 
we obtain that $f$ is satisfiable if and only if 
the modified graph $G'$ satisfies $\rho_{k_1}(G')=\rho_{k_2}(G')$.

Note that we require $\lceil(3k_2+1)/2\rceil\leq 2k_1+1$
instead of just $k_2<2k_1$ in order to ensure that
$\rho_{k_1}(G(C_j))=2$.
$\Box$

\medskip

\noindent We proceed to consequences of Theorem \ref{theorem1} for distance matching numbers. Note that a graph $G$ satisfies $\nu_k(G)=\nu_{2k}(G)$  if and only if 
$\rho_1(L(G)^k)=\rho_2(L(G)^k)$, 
that is, these graphs can be recognized by a very simple algorithm.

For a positive integer $k$, 
a {\it $k$-unit} is a pair $(G,e)$ 
where $G$ is a connected graph, $e$ is an edge of $G$, and $\nu_k(G)=1$.
The {\it boundary $\partial (G,e)$ of $(G,e)$}
is the set of vertices of $G$ that are at distance exactly $k$ from $e$ in $G$.
Note that, since $G$ is connected and $\nu_k(G)=1$, no vertex of $G$ is at distance more than $k$ from $e$ in $G$, 
and the boundary $\partial (G,e)$ is independent.

\begin{corollary}\label{corollary2}
A graph $G$ satisfies $\nu_k(G)=\nu_{2k}(G)$ for some positive integer $k$ if and only if 
$G$ arises from the disjoint union of $k$-units
$(G_1,e_1),\ldots,(G_\ell,e_\ell)$ by arbitrarily identifying vertices 
in $\bigcup_{i=1}^\ell\partial (G_i,e_i)$, where $\ell=\nu_k(G)$.
\end{corollary}
{\it Proof:}  Let $G$ be a graph.

In order to prove the sufficiency, let $G$ arise in the described way from the $k$-units $(G_i,e_i)$.
Let $P$ be a maximum $1$-packing of $L(G)^k$, 
that is, $P$ is a set of edges of $G$ that are at pairwise distance more than $k$ in $L(G)$.
Since $\nu_k(G_i)=1$, the set $P$ contains at most one edge from each $G_i$,
which implies that $\nu_k(G)=\rho_1(L(G)^k)=|P|\leq \ell$.
By the definition of the boundary, the set $\{ e_i:i\in [\ell]\}$ is a $2$-packing of $L(G)^k$,
and hence $\nu_{2k}(G)=\rho_2(L(G)^k)\geq \ell$, which implies $\nu_k(G)=\nu_{2k}(G)$.

In order to prove the necessity, let $G$ satisfy $\nu_k(G)=\nu_{2k}(G)$,
that is, $\rho_1(L(G)^k)=\rho_2(L(G)^k)$.
Let $P$ be a maximum $2$-packing of $L(G)^k$.
By Theorem \ref{theorem1}, the set $P$ is a transversal of ${\cal S}(L(G)^k)$
and $\{ N_{L(G)^k}[e]:e\in P\}$ is a partition of $E(G)$, the vertex set of $L(G)^k$, 
into sets that are complete in $L(G)^k$.
Let $P=\{ e_1,\ldots,e_\ell\}$ and let $E_i=N_{L(G)^k}[e_i]$ for $i\in [\ell]$.
For $i\in [\ell]$, 
let $V_i$ denote the set of vertices of $G$ that are incident with an edge in $E_i$,
and let $G_i=(V_i,E_i)$.
By definition, and since $E_i$ is a complete set in $L(G)^k$,
the graph $G_i$ is connected, $e_i$ is an edge of $G_i$, and  $\nu_k(G_i)=1$,
that is, $(G_i,e_i)$ is a $k$-unit.
Note that the graphs $G_i$ are edge-disjoint yet not vertex-disjoint subgraphs of $G$.
If $G_i$ and $G_j$ share a vertex $u$ for some $i\not=j$, and $u$ does not belong to the intersection of the boundaries
$\partial (G_i,e_i)\cap \partial (G_j,e_j)$,
then the distance in $L(G)^k$ between $e_i$ and $e_j$ is at most $2$,
which is a contradiction.
Hence $G$ arises in the described way 
from the $k$-units $(G_i,e_i)$. 
This completes the proof.
$\Box$

\medskip

\noindent Let $G$ be a graph.
A vertex of degree $1$ in $G$ is a {\it leaf of $G$}.
A triangle $uvwu$ in $G$ such that the degree of $u$ and $v$ in $G$ is $2$
is a {\it pendant triangle of $G$} and the edge $uv$ is a {\it triangle edge of $G$}.

\begin{corollary}[Cameron and Walker \cite{cawa}]\label{corollary3}
A connected graph $G$ satisfies $\nu_1(G)=\nu_2(G)$ if and only if $G$ is either 
a star, or a triangle, or
arises from a connected bipartite graph with two non-empty partite sets $V_1$ and $V_2$ by 
\begin{itemize}
\item attaching at least one and possibly more leaves to each vertex in $V_1$, and
\item attaching pendant triangles to some vertices in $V_2$.
\end{itemize}
\end{corollary}
{\it Proof:}  Let $G$ be a graph that satisfies satisfies $\nu_1(G)=\nu_2(G)$.
By Corollary \ref{corollary2},
the graph $G$ arises from the disjoint union of $1$-units 
by arbitrarily identifying vertices in their boundaries.
It follows immediately from the definition 
that if $(G,e)$ is a $1$-unit, then 
\begin{itemize}
\item either $G$ is a star and $\partial (G,e)$ is the set of leaves of $G$
that are not incident with $e$,
\item or $G$ is a triangle and $\partial (G,e)$ consists of the vertex 
that is not incident with $e$.
\end{itemize}
The desired structure not follows immediately.
In fact, $V_1$ is the set of all centers of $1$-units that are stars 
and $V_2$ is the union of all boundaries (after identification). $\Box$

\medskip

\noindent Our results motivate some questions.
In view of Observation \ref{observation1} it might make sense to consider bounds for $\frac{\rho_{k_1}(G)}{\rho_{k_2}(G)}$ rather than linear relations between $\rho_{k_1}(G)$ and $\rho_{k_2}(G)$.
It would be interesting to know whether the decision problems considered in Theorem \ref{theorem2} and Corollary \ref{coro NPhard} are in NP. 
We believe that 
for all positive integers $k_1$ and $k_2$ with $k_1<k_2<2k_1$,
it is NP-hard to determine for a given graph $G$
whether $\rho_{k_1}(G)=\rho_{k_2}(G)$.
Unfortunately, Corollary \ref{coro NPhard} does not cover all possible cases.

\end{document}